\documentclass[11pt]{amsart}
\usepackage{amssymb}

\theoremstyle{plain}
\newtheorem{theorem}{Theorem }[section]

\newtheorem{lemma}[theorem]{Lemma}
\newtheorem{corol}[theorem]{Corollary}
\newtheorem{defi}[theorem]{Definition}
\theoremstyle{remark}
\newtheorem{rem}{Remark}

\newcommand{\mundo}{\Omega \times M}
\newcommand{\pot}{L^1(\Omega,C(M))}
\newcommand{\SA}{\mathcal{A}}
\newcommand{\SB}{\mathcal{B}}
\newcommand{\SC}{\mathcal{C}}
\newcommand{\SD}{\mathcal{D}}
\newcommand{\SM}{\mathcal{M}}
\newcommand{\SF}{\mathcal{F}}
\newcommand{\SK}{\mathcal{K}}

\newcommand{\SME}{\mathcal{M}_e}
\newcommand{\SP}{\mathcal{P}}

\newcommand{\la}{\lambda}
\newcommand{\Om}{\Omega}
\newcommand{\BP}{\mathbb{P}}
\newcommand{\R}{\mathbb{R}}

\newcommand{\N}{\mathbb{N}}

\newcommand{\eps}{\varepsilon}

\newcommand{\bt}{\bigtriangleup}
\newcommand{\real}{\mathcal{R}}

\begin{document}

\title[Equilibrium States for Random Non-uniformly Expanding Maps]{Equilibrium
States for Random Non-uniformly Expanding Maps}
\date{March 12, 2003}
\author[A. Arbieto, C. Matheus and K. Oliveira]{Alexander Arbieto, Carlos Matheus and Krerley
Oliveira}
\thanks{A.A. was partially supported by CNPq-Brazil, C.M. was partially supported by
Faperj-Brazil and K.O. was partially supported by CNPq-Brazil}

\begin{abstract}We show that, for a robust ($C^2$-open) class of random
non-uniformly expanding maps, there exists equilibrium states for a large class
of 
potentials.In particular, these sytems have measures of maximal entropy. These
results also give a partial answer to a question posed by Liu-Zhao. The proof of the main
result uses an extension of techniques in recent works by Alves-Ara\'ujo, 
Alves-Bonatti-Viana and Oliveira. 
 
\end{abstract}

\maketitle

\section{Introduction}
\label{intro}

Particles systems, as they appear in kinetic theory of gases, have been an
important model motivating much development in the field of Dynamical Sytems and
Ergodic Theory. While these are deterministic systems, ruled by Hamiltonian
dynamics, the evolution law is too complicated, given the huge number of
particles involved. Instead, one uses a stochastic approach to such systems.

More generally, ideais from statistical mechanics have been brought to the
setting of dynamical systems, both discrete-time and continuous-time, by Sinai,
Ruelle, Bowen, leading to a beautiful and very complete theory of equilibrium
states for uniformly hyperbolic diffeomorphisms and flows. In a few words,
equilibrium states are invariant probabilities in the phase space which
maximaze a certain variational principle (corresponding to the Gibbs free energy
in the statistical mechanics context). The theory of Sinai-Ruelle-Bowen gives
that for uniformly hyperbolic systems equilibrium states exist, and they are
unique if the system is topologically
transitive and the potential is H\"older continuous.

Several authors have worked on extending this theory beyond the uniformly
hyperbolic case. See e.g.~\cite{Buzzi},~\cite{Sarig}, among other important
authors.
Our present work is more directly motivated by the results of Oliveira~\cite{O}
where he constructed equilibrium states associated to potentials with
not-too-large variation, for a robust ($C^1$-open) class of non-uniformly
expanding maps introduced by Alves-Bonatti-Viana~\cite{ABV}.

On the other hand, corresponding problems have been studied also in the context
of the theory of random maps, which was much developed by
Kifer~\cite{K} and Arnold~\cite{Arnold}, among other mathematicians.
Indeed, Kifer~\cite{K} proved the existence of equilibrium states for random
uniformly exapnding systems, and Liu~\cite{L} extended this to uniformly
hyperbolic systems. 

In the present work, we combine these two approaches to give a construction of
equilibrium states for non-uniformly hyperbolic maps. In fact, some attempts to
show the existence of equilibrium states beyond uniform hyperbolicity were made by Khanin-Kifer~\cite{KK}.
However, our point of view is quite different.
Before stating the main result, we
recall that a random map is a continuous map $f:\Omega\rightarrow C^r(M,M)$
where $M$ is a compact manifold $\Omega$ is a Polish space, and
$T:\Omega\rightarrow\Omega$ a measurably invertible continuous map with an
invariant ergodic measure $\BP$. 
The main result is the following :

\begin{themaintheorem} {\it``For a $C^2$-open set $\SF$ of non-uniformly expanding local 
diffeomorphisms, potentials $\phi$ with \emph{low variation} and
$f:\Omega\rightarrow\SF$, 
there are equilibrium states for the random system associated to $f$ and $T$. 
In particular, $f$ admits measures with maximal entropy.''}
\end{themaintheorem}      

A potential has low variation if it is not far from being constant. See
the precise definition in section $3$. In particular, constant functions have
low variation; their equilibrium states are measures of maximal entropy. 

The proof, which we present in the next sections extends ideias from
Alves-Ara\'ujo~\cite{AA}, Alves-Bonatti-Viana~\cite{ABV} and Oliveira~\cite{O}.

It is very natural to ask whether these equilibrium states we construct are
unique and whether they are (weak) Gibbs states, Another very interesting
question is whether existence (and uniqueness) of equilibrium states extends to
(random r deterministic) non-uniformly hyperbolic maps with singularities, such
as the Viana maps~\cite{AA}. Although our present methods do not solve these
questions, we believe the answers are affirmative.

\section{Definitions}\label{s.12}

{\it Random Transformations and Invariant Measures}

\smallskip

Let $M^l$ be a compact $l$-dimensional Riemannian manifold and $\SD$ the space of 
$C^{2}$ local diffeomorphisms of $M$. 
Let $(\Om ,T , \BP)$ a measure preserving system,  where $T : \Om \rightarrow
\Om$ is $\BP$-invariant ($\BP$ is a Borel measure) and $\Omega$ is a Polish space, i.e., $\Omega$ is a
complete separable metric space. 
By a random transformation we understand a continuous map 
$f:\Om \rightarrow \SD$. Then we define:
\begin{equation}\label{e.cociclo}
f^n(w)=f(T^{n-1}(w))\circ \cdots \circ f(w), \textrm{      }f^{-n}(w)=(f^n(w))^{-1}.
\end{equation}

We also define the skew-product generating by $f$:
$$F:\mundo  \rightarrow \mundo, \textrm {   } F(w,x)=(Tw, f(w)x).$$

We denote $\SP(\mundo)$ the space of probability measures $\mu$ on $\mundo$
such that the marginal of $\mu$ on $\Om$ is $\BP$. 
Let $\SM(\mundo)\subset \SP(\mundo)$ be the measures $\mu$ which are $F$-invariant.

Because $M$ is compact, invariant measures always exists and the property of $\BP$ be the marginal on $\Om$ of a invariant measures can be characterized by its disintegration:
$$d\mu(w,x)=d\mu_w(x)d\BP(w).$$
$\mu_w$ are called {\it samples measures} of $\mu$ (see~\cite{Liu},~\cite{LQ}).

An invariant measure is called {\it ergodic} if $(F,\mu)$ is ergodic, the set of all ergodic measures is denoted by $\SME(\mundo)$. Furthermore, each invariant measure can be decomposed into its ergodic components by integration when the $\sigma$-algebra on $\Om$ is countably generated and $\BP$ is ergodic.                                                                                                                                                                                                                                                                                                                                                                                                                                                                                                           

\textbf{In what follows}, as usual, \textbf{we always assume $(\Omega
,\SA,\BP)$ is a Lebesgue space, $(T,\BP)$
is ergodic and $T$ is measurably invertible and continuous}. Observe that these 
assumptions
are satified in the canonical case of left-shift operators $\tau$, $\Omega$
being $C^r(M,M)^{\mathbb{N}}$ or $C^r(M,M)^{\mathbb{Z}}$.

\smallskip

{\it Entropy}

\smallskip

We follow Liu~\cite{Liu} on the definition of the Kolmogorov-Sinai entropy for random transformations:

Let $\mu$ an $F$-invariant measure like above. Let $\xi$ a finite Borel partition of $M$. We set:

\begin{equation}
\label{e.entropy}
h_{\mu}(
f,\xi)=\lim_{n\to +\infty}\frac{1}{n}\int
H_{\mu_w}(\vee_{k=0}^{n-1}f^k(w)^{-1}\xi)d\BP(w),
\end{equation}

where $H_{\nu}(\eta):=-\sum_{C\in \eta}\nu(C)\log \nu(C)$ (and $0\log 0=0$), for
a finite partition $\eta$ and $\nu$ a probability on $M$ (and $\mu_w$ are the
sample measures of $\mu$).

\begin{defi}
The entropy of $(f,\mu)$ is: 
$$h_{\mu}(f):=\sup_{\xi}h_{\mu}(f,\xi)$$with the supremum taken over all finite Borel partitions of $M$.
\end{defi}

\begin{defi} The \emph{topological entropy} of $f$ is
$h_{top}(f)=\sup\limits_{\mu}h_{\mu}(f)$ 
\end{defi}

\begin{theorem}[``Random'' Kolmogorov-Sinai theorem]\label{Kolmogorov-Sinai} If
$\SB$ is the Borel $\sigma$-algebra of $M$ and $\xi$ is a generating partition
of $M$,
i.e.,
$$
\bigvee_{k=0}^{+\infty} f^{-k}(w) \ \xi = \SB  \ \textrm{ for }\BP-a.e. \ w
,$$
then
$$
h_\mu(f)=h_\mu (f, \xi ) .
$$

\end{theorem}

For a proof of this theorem see~\cite{LQ} or~\cite{B}.

\newpage

\smallskip

{\it Equilibrium States}

\smallskip

Let $\pot$ the set of all families $\{\phi=\{\phi_w\in C^0(M)\}\}$ such that the map 
$(w,x)\rightarrow \phi_w(x)$ is a measurable map and $\|\phi\|_1:=\int_{\Om} |\phi_w|_{\infty}d\BP(w) <+\infty$.

For a $\phi \in \pot$, $\eps>0$ and $n\geq 1$, we define:

$$\pi_f(\phi)(w,n,\eps)=\sup \{ \sum_{x\in K} e^{S_f(\phi)(w,n,x)};  K \textrm{ is a }
(n,\eps)-\textrm{separated set}\},$$
where $S_f(\phi)(w,n,x):=\sum\limits_{k=0}^{n-1}\phi_{T^k(w)}(f^k(w)x)$.

\begin{defi}
The map $\pi_f:\pot \rightarrow \R\cup \{\infty\}$ given by:
$$\pi_f(\phi)=\lim_{\eps\to 0} \limsup_{n\to +\infty} \frac{1}{n}\int_{\Om}\log
\pi_f(\phi)(w,n,\eps) d\BP(w).$$
is called {\it the pressure map}. 
\end{defi}

It is well know that the variational principle occurs (see~\cite{Liu}):

\begin{theorem}
\label{t.var}
If $\Om$ is a Lebesgue space, then for any $\phi \in \pot$ we have:
\begin{equation}
\label{e.pressure}
\pi_f(\phi)=\sup_{\mu\in\SM(\mundo)} \{h_{\mu}(f)+\int\phi d\mu\}
\end{equation}
\end{theorem}

\begin{rem}\label{r.erg}
If $\BP$ is ergodic then we can take the supremum over the set of ergodic measures.
\end{rem}

\begin{defi}
A measure $\mu \in \SM(\mundo)$ is an equilibrium state for $
f$, if $\mu$ attains the supremum of~(\ref{e.pressure}). 
\end{defi}

\smallskip

{\it Physical Measures}

\smallskip

As in the deterministic case, we follow~\cite{AA} on the definition of
\emph{physical measure} in the context of random transformations :

\begin{defi} A measure $\mu$ is a physical measure if for positive Lebesgue
measure set of points $x\in M$ (called the \emph{basin} $B(\mu)$ of $\mu$),

\begin{equation}
\label{medida fisica}
\lim\limits_{n\rightarrow\infty}\frac{1}{n}\sum_{j=1}^{n-1}\phi(f^j(w)(x)) =
\int\phi d\mu \textrm{ for all continuous } \phi :M\rightarrow\mathbb{R} .
\end{equation}

for $\BP$-ae $w$.
\end{defi}

\smallskip

\section{Statement of the results} \label{s.2}

Before starting abstract definitions, we comment that in next section, it is showed
that there are examples of random transformations satisfying our hypothesis below. 

We say that a local diffeomorphism $f$ of $M$ is in $\widetilde{\SF}$ if $f$ is in $\SD$ and satisfies, for positive constants $\delta_0$, $\beta$, $\delta_1$, $\sigma_1$, and $p$, $q\in\N$,
the following properties :

\begin{enumerate}

\item[(H1)] There exists a covering $B_1, \dots,B_p, \dots, B_{p+q}$ of
      $M$ such that every $f|B_i$ is injective and 

\begin{itemize}

\item  $f$ is uniformly expanding at every $x\in B_{1} \cup \dots \cup
B_{p}$: $$\|Df(x)^{-1}\| \leq (1+\delta_1)^{-1}.$$

\item  $f$ is never too contracting: $\|Df(x)^{-1}\|\leq (1+\delta_0)$ for
every $x\in
M$.

\end{itemize}

\item[(H2)] $f$ is everywhere volume-expanding:  $|\det Df(x)|\geq \sigma_1$
      with $\sigma_1>q$.

Define $$V=\{x\in M; \|Df(x)^{-1}\|>(1+\delta_1)^{-1}\}.$$

\item[(H3)] There exists a set $W \subset B_{p+1}\cup\dots\cup B_{p+q}$ containing $V$ such that 
$$
M_1>m_2 \quad\text{and}\quad m_2-m_1<\beta
$$
where $m_1$ and $m_2$ are the infimum and the supremum 
of $|\det Df|$ on $V$, respectively, 
and $M_1$ and $M_2$ are the
infimum and the supremum of $|\det Df|$ on $W^c$, respectively. 
\end{enumerate}

This kind of transformations was considered by~\cite{ABV}, ~\cite{O}, ~\cite{AA}, where they construct $C^1$-open sets of such maps.

We will consider a subset $\SF\subset\widetilde{\SF}$ such that :

\begin{enumerate}

\item[(C1)] There is a uniform constant $A_0$ s.t.
$|\log \| f\|_{C^2}|\leq A_0$ for any $f\in\SF$ and the constants $m_1,m_2,
M_1, M_2$ are uniform on $\SF$ ;

\end{enumerate}

\textbf{From now on, our random transformations will be given by a map $f:\Om \rightarrow
\SF$}, and $f$ satisfy the following condition :

\begin{enumerate}

\item[(C2)] $f$ admits an ergodic absolutely
continuous \emph{physical
measure} $\mu_{\BP}$ (see section 2).

\end{enumerate} 

\begin{rem} We will show in the appendix that $(H1), (H2)$ implies the
following property:

\begin{enumerate}

\item[(F1)] There exists some $\gamma_0= \gamma_0(\delta_1, \sigma_1,
p, q)<1$ such that the random orbits of Lebesgue almost every point spends at most a fraction
of time $\gamma_0 <1$ inside $B_{p+1}\cup\dots\cup B_{p+q}$, 
depending only on $\sigma_1$, $p$, $q$. I.e., for $\BP$-a.e. $w$ and Lebesgue almost
every $x$
$$\lim\limits_{n\rightarrow\infty}\frac{ \# \{ 0\leq j\leq n-1: f^j(w)(x)\in
B_{p+1}\cup\dots\cup B_{p+q}\} }{n}\leq\gamma_0 .$$

\end{enumerate}

\end{rem}

Then we analyse the existence of an equilibrium state for low-variation potentials:

\begin{defi}
\label{d.potencial}
A potential $\phi \in \pot$ has $\rho_0$-low variation if

\begin{equation}
\label{e.low}
\| \phi\|_1 <
\pi_{f}(\phi)-\rho_0 h_{top}(f).
\end{equation}
\end{defi}

\begin{rem} We call $\phi$ above a $\rho_0$-low variation potential because in
the deterministic case (i.e., $\phi(w,x)=\phi(x)$), if $\max\phi -\min\phi <
(1-\rho_0)h_{top}(f)$ then $\phi$ satisfies~(\ref{e.low}).
\end{rem}

The main result is :
\bigskip

\textbf{Theorem A.}
\textit{Assume hypotheses (H1), (H2), (H3) hold, with $\delta_0$ and
$\beta$ sufficiently small and assume also conditions (C1), (C2). Then, there exists $\rho_0$ such that if
$\phi$ is a continuous potential with $\rho_0$-low variation then
$\phi$ has some equilibrium state. Moreover, these equilibrium
states are hyperbolic measures, with all Lyapunov exponents bigger
than some $c=c(\delta_1, \sigma_1, p, q)>0$. }


\section{Examples}
\label{exemplos}

In this section we exhibit a $C^1$-open class of $C^2$-diffeomorphism which are
contained in $\widetilde{\SF}$. To start the construction, we now follow~\cite{O}
\emph{ipsis-literis} and construct examples of `deterministic' non-uniformly
expanding maps. After this, we construct the desired random non-uniformly
expanding maps in $\SF$ a $C^2$-neighborhood of a fixed diffeomorphism of
$\widetilde{\SF}$.

We observe that the class $\SF$ contains an open set of non-uniformly exapanding
which \emph{are not uniformly exapnding}.

We start by considering any Riemannian manifold that supports an expanding map $g
:M \rightarrow M$. For simplicity, choose  $M=\mathbb{T}^n$
  the $n$-dimensional torus, and $g$ an endomorphism induced from a linear
map with eigenvalues $\lambda_n>\dots>\lambda_1>1.$ Denote by $E_i(x)$ the eigenspace
associated to the eigenvalue $\lambda_i$ in $T_x M$.  
 
Since $g$ is an expanding map,   $g$ admits a transitive Markov partition $ R_1, \dots, R_d $ with arbitrary small
  diameter. We may suppose that $g|R_i$ is injective for every $i=1,\dots,d$. Replacing   by a iterate if necessary,  we may suppose that there exists a fixed
point $p_0$ of $g$ and, renumbering if necessary,  this point is contained in the interior of the rectangle $R_d$ of the Markov partition.

Considering a small neighborhood $W \subset R_d$ of $p_0$ we 
deform $g$ inside $W$ along the  direction $E_1$. This deformation consists essentially in rescaling the expansion
along the invariant manifold associated to $E_1$ by a real function $\alpha$.  Let us be more precise:

Considering $W$ small, we may identify $W$ with a neighborhood of
$0$ in $\real^n$ and $p_0$ with $0$. Without loss of generality, suppose that $W=(-2\epsilon,2\epsilon)
\times B_{3r}(0)$, 
where $B_{3r}(0)$ is the ball or radius $3r$ and
center $0$ in $\real^{n-1}$. Consider a function $\alpha:(-2\epsilon,2\epsilon)\rightarrow  \real$ such
$\alpha(x)=\lambda_1x$ for every $|x|\geq \epsilon$ and for small constants $\gamma_1,\gamma_2$:

\begin{enumerate}

\item  $(1+\gamma_1)^{-1}<\alpha'(x)< \lambda_1 + \gamma_2$ 

\item $\alpha'(x)<1$ for every $x \in (-\frac{\epsilon}{2},\frac{\epsilon}{2})$; 

\item $\alpha$ is $C^0$-close to $\lambda_1$: $\sup\limits_{x\in (-\epsilon,\epsilon)}|\alpha(x)-\lambda_1x| < \gamma_2$,

\end{enumerate}

Also, we consider a bump function $\theta: B_{3r}(0)\rightarrow \real$ such $\theta(x)=0$ for every $2r\leq |x| \leq 3r$
 and $\theta(x)=1$ for every $0\leq |x|\leq r$. Suppose that $\|\theta'(x)\|\leq C$ for
every $x\in B_{3r}(0)$.  Considering coordinates $(x_1,\dots,x_n)$ such that $\partial_{x_i} \in E_i$, define
$f_0$ by:

$$ f_0(x_1,\dots,x_n)=
(\lambda_1x_1+\theta(x_2,\dots,x_n)(\alpha(x_1)-\lambda_1x_1), \lambda_2x_2,\dots,\lambda_n x_n)$$

Observe that by the definition of $\theta$ and $\alpha$ we can extend $f_0$
smoothly to $\mathbb{T}^n$ as $f_0=g$
outside $W$. Now, is not difficult to prove that $f_0$ satisfies the conditions (H1), (H2), (H3) above. 

First, we have that $\|Df_0(x)^{-1}\|^{-1}\geq \min\limits_{i=1,\dots,n} \|\partial_{x_i} f_0\|.$ Observe that:

$$\partial_{x_1} f_0(x_1,\dots,x_n)= (\alpha'(x_1)\theta(x_2,\dots,x_n)+(1-\theta(x_2,\dots,x_n))\lambda_1,0,\dots,0)$$  
$$\partial_{x_i} f_0(x_1,\dots,x_n) = ((\alpha(x_1)-\lambda_1)\partial_{x_i}\theta(x_2,\dots,x_n),0,\dots,\lambda_i,0,\dots,0), \text{ for } i\geq 2.$$

Then, since $\|\partial_{x_i} \theta(x)\| \leq C$ for every $x\in B_{3r}(0)$,  and $\alpha(x_1)-\lambda_1x_1
\leq \gamma_2$ we have that $ \|\partial_{x_i} f_0\|>(\lambda_i-\gamma_2C)$ for every
$i=2,\dots,n$. Moreover, by condition 1, $\|\partial_{x_1} f_0\| \leq \max \{\alpha'(x_1), \lambda_1\}\leq
\lambda_1 + \gamma_2,$ if we choose $\gamma_2$ small in such way that $\lambda_2-\gamma_2C >
\lambda_1 +\gamma_2 $ then:

$$\|\partial_{x_i} f_0\|>\|\partial_{x_1} f_0\|, \text{ for every } i\geq 2.$$ 

Notice also that $\|\partial_{x_1} f_0\| \geq \min\{\alpha'(x_1),\lambda_1\} \geq (1+\gamma_1)^{-1}.$ This prove
that: 

$$\|Df_0(x)^{-1}\|^{-1}\geq \min\limits_{i=1,\dots,n} \|\partial_{x_i} f_0\| \
(1+\gamma_1)^{-1}.$$ Since $f$ coincides with $g$ outside $W$, we have  $\|Df_0(x)^{-1}\|\leq \lambda_1^{-1}$ for every
$x\in W^c$. Together with the above inequality, this proves  condition (H1), with $\delta_0=\gamma_1$.    
     
Choosing $\gamma_1$ small and $p=d-1$, $q=1,$ $B_i=R_i$ for every $i=1,\dots,d$,  condition (H2) is immediate.
Indeed, observe that the Jacobian of $f_0$  is given by the formula:

 $$\det Df_0(x)= (\alpha'(x_1)\theta(x_2,\dots,x_n)+(1-\theta(x_2,\dots,x_n))\lambda_1)\prod_{i=2}^n
 \lambda_i.$$  Then,  if we choose $\gamma_1 < \prod_{i=2}^n
 \lambda_i-1$:
 
 $$\det Df_0(x)>(1+\gamma_1)^{-1}\prod_{i=2}^n
 \lambda_i >1.$$ Therefore, we may take $\sigma_1= (1+\gamma_1)^{-1}\prod_{i=2}^n \lambda_i>1.$

To verify property (H3) for $f_0$, observe that if we denote by 

$$V=\{x\in M; \|Df_0(x)^{-1}\|>(1+\delta_1)^{-1}\},$$ 

with $\delta_1<\lambda_1-1,$ then $V\subset W.$ Indeed, since $\alpha(x_1)$ is constant equal to
$\lambda_1x_1$ outside $W$  we have that $\|Df_0(x)^{-1}\|\leq \lambda_1^{-1} < (1+\delta_1)^{-1}$, for every
$x\in W^c$.  Given $\gamma_3$ close to 0, we may choose $\delta_1$ close to 0 and $\alpha$ satisfying the conditions above in 
such way that,

$$\sup\limits_{x,y\in V}\alpha'(x_1)-\alpha'(y_1)< \gamma_3.$$ If $m_1$ and $m_2$ are the 
infimum and the supremum of $|\det Df_0|$ on $V$, respectively, 

$$m_2-m_1\leq C(\sup\limits_{x,y\in V}
\alpha'(x_1)-\alpha'(y_1)) < \gamma_3 C,$$
where $C=\prod\limits_{i=2}^n\lambda_i$. 
Then, we may take $\beta=\gamma_3C$ in (H3). If $M_1$ is the 
infimum  of $|\det Df_0|$ on $W^c$, $M_1>m_2$, since $\lambda_1> (1+\delta_1)\geq \sup\limits_{x\in V}
\alpha'(x)$. 

The arguments above show that the hypotheses $(H1), (H2),(H3)$ are satisfied by $f_0$. Moreover, if we one takes
    $\alpha(0)=0$, then $p_0$ is fixed point for $f_0$, which is not a repeller, since $\alpha'(0)<1$.
Therefore, $f_0$ is not a uniformly expanding  map.  
    
It is not difficult to see that this construction may be carried out in such way that $f_0$ does not satisfy the  expansiveness
 property: there is a fixed hyperbolic saddle point $p_0$ such that the stable manifold of $p_0$ is 
 contained in the unstable manifold of two other fixed points.

Now, if $\SF$ denotes a small $C^2$-neighborhood of $f_0$ in $\widetilde{\SF}$, and
$h:\Omega\rightarrow\widetilde{\SF}$ is a continuous map, Alves-Ara\'ujo~\cite{AA} shows
that if $w^*\in\Omega$ is such that $f(w^*)=f_0$ and $\theta_{\eps}$ is a
sequence of measures, $\textrm{supp}(\theta_{\eps})\rightarrow \{ w_0\}$ then for
small $\eps>0$ there are physical measures for the RDS $f:\Omega^{\mathbb{Z}}
\rightarrow\SF$, $f(\dots,w_{-k},\dots, w_0,\dots ,w_k,\dots)= h(w_0)$ . This
concludes the construction of examples satisfying $(H1),(H2),(H3),(C1),(C2)$.


\section{Proof of the theorem A}
\label{prova}

We now precise the conditions on $\delta_0$ and $\beta$. We consider $\gamma_0$
given in condition $(F1)$. By condition $(C1)$, there exists
$\epsilon_0 >0$ s.t. for any $\eta\in B_{\epsilon_0}(\xi)$ and
$\BP$-a.e. $w$ holds : 
$$\frac{|| Df(w)^{-1}(\xi) ||}{|| Df(w)^{-1}(\eta) ||}\leq e^{\frac{c}{2}},$$
where $c$ is such that for some $\alpha>\gamma_0$, we have $(1+\delta_0)^{\alpha}(1+\delta_1)^{-(1-\alpha )}<
e^{-2c}<1$ and $\alpha m_2 + (1-\alpha )M_2<\gamma_0 m_1+ (1-\gamma_0)M_1
-l\log(1+\delta_0)$ ($l:=\dim (M)$), if $\delta_0$ and $\beta$ are sufficiently
small. Now, the constants fixed above allows us to
prove good properties for the objects defined below, which are of fundamental interest in the proof of
theorem A.
 
\smallskip
{\it Expansive Measures and Hyperbolic Times}

\begin{defi}
We say that a measure $\nu \in \SM({\mundo})$ is \emph{$f$-expanding with
exponent $c$} if for $\nu$-almost every $(w,x)\in \mundo$ we have:
$$
\la(w,x)=\limsup\limits_{n\rightarrow +\infty}
\frac{1}{n}\sum\limits_{j=0}^{n-1} \log \|Df(T^j(w))(f^j(w)(x))^{-1}\| \leq
-2c <0.
$$
\end{defi}

\begin{defi}
We say that $n$ is a \textit{hyperbolic time} for $(w,x)$ with exponent $c$, if
for
every  $1\leq k\leq n$:
$$
\prod _{j=n-k}^{n-1}  \|Df(T^{j+1}(w))(f^j(w)(x))^{-1}\| \leq e^{-ck}.
$$
\end{defi}

As in lemma 3.1 of~\cite{ABV}, lemma 4.8 of~\cite{O} and lemma 2.2 of~\cite{AA},
we have infinity many hyperbolic times for expanding measures. For this we need
a lemma due to Pliss (see~\cite{ABV}).

\begin{lemma}
\label{l.pliss}

Let $A\geq c_2>c_1>0$ and $\zeta=\frac{c_2-c_1}{A-c_1}$. Given real numbers $a_1,\cdots,a_N$ satisfying:
$$\sum_{j=1}^{N}a_j \geq c_2N \textrm{  and  } a_j\leq H \textrm{ for all } 1\leq j\leq,$$ there are $l>\zeta N$ and $1<n_1<\cdots<n_l\leq N$ such that:
$$\sum_{j=n+1}^{n_i}a_j \geq c_1(n_1-n) \textrm{ for each } 0\leq n<n_i,\textrm{ }i=1,\cdots ,l.$$
\end{lemma}

\begin{lemma}
\label{l.hyptimes}
For every invariant measure $\nu$
which is $c$-expanding, there exists a full
$\nu$-measure set $H\subset \Omega\times M$ such that every $(w,x)\in H$ has
infinitely many hyperbolic times $n_i=n_i(w,x)$ with exponent $c$ and, in
fact, the
density of hyperbolic times at infinity is larger than some
$d_0=d_0(c)>0$:

\begin{enumerate}
\item $\displaystyle\prod _{j=n-k}^{n-1}  \|Df(T^{j+1}(w))(f^j(w)(x))^{-1}\|
\leq e^{-cj}$ for every $1\leq k\leq n_i$
\item $\displaystyle{\liminf\limits_{n\rightarrow \infty}
\frac{\sharp\{0\leq n_i\leq n\}}{n}\geq d_0>0}$.
\end{enumerate}
\end{lemma}

\begin{proof}
Let $H\subset \mundo$ with full $\nu$-measure. For any $(w,x)\in H$ and $n$
large enough, we have:
$$\sum\limits_{j=0}^{n-1} \log \|Df(T^j(w))(f^j(w)(x))^{-1}\| \leq
-\frac{3c}{2}n  $$
Now, by $(C1)$ we can apply lemma~\ref{l.pliss} with $A = \sup\limits_{(w,x)}
(-\log ||Df(w)^{-1}(x)||)$,
$c_1=c$, $c_2=\frac{3c}{2}$ and $a_i=-\log \|Df(T^j(w))(f^j(w)(x))^{-1}\|$ and
the statement follows.
\end{proof}

\begin{lemma}\label{l.var} $\exists \ \epsilon_0 >0$
such that for $\BP$-a.e.~$w$, if $n_i$ is a hyperbolic time of $(w,x)$ 
and $f^{n_i}(w)(z)\in B_{\epsilon_0}(f^{n_i}(w)(x))$ then $d(f^{n_i -j}(w)(z),
f^{n_i -j}(w)(x))\leq e^{\frac{-cj}{2}} d(f^{n_i}(w)(z),f^{n_i}(w)(x))$, $\forall
1\leq j\leq n_i$.
\end{lemma}

\begin{proof} By $(C1)$ we know that there exists $\eps_0>0$ such that for any $\eta\in
B_{\eps_0}(\xi)$ we have:
$$\frac{|| Df(w)^{-1}(\xi) ||}{|| Df(w)^{-1}(\eta) ||}\leq e^{\frac{c}{2}}
\textrm{  for $\BP$-ae $w$}.$$
In fact, this hold in the $T$-orbit of $w$ $\BP$-ae. Indeed, let $C=\{w;$ 
$\frac{|| Df(w)^{-1}(\xi) ||}{|| Df(w)^{-1}(\eta) ||}\leq e^{\frac{c}{2}}\}$ for
any $\xi$ and $\eta\in B_{\eps_0}(\xi)$, then $\bigcap T^j(C)$ has full measure
and the estimate follows. Because $f^{n_i}(w)(z)\in
B_{\epsilon_0}(f^{n_i}(w)(x))$, by the estimative above, we have that $w$ $\BP$-ae if we take the
inverse branch of $f^{n_i}(w)$ which sends $f^{n_i}(w)(x)$ to $f^{n_i-1}(w)(x)$
(restricted to $B_{\eps_0}(f^{n_i}(w)(x)$)) and has derivative with norm less than
$e^{-\frac{c}{2}}$, then we have $d(f^{n_i-1}(w)(z), f^{n_i-1}(w)(x) ) \leq
\eps_0$. Using the estimate along the orbit (and induction), we have:
$$\prod _{j=n-k}^{n-1}  \|Df(T^{j+1}(w))(f^j(w)(z))^{-1}\|
\leq e^{-\frac{ck}{2}}\textrm{ for all $0\leq k \leq n_i$.}$$
The statement follows.
\end{proof}

Now we define a set of measures where the ``bad set'' $V$ has small measure.

\begin{defi} We define the convex set $K_{\alpha}$ by 
\begin{center}
$K_{\alpha}=\{ \mu : \mu (\Omega\times V)\leq\alpha  \}$
\end{center}
\end{defi}

\begin{lemma}\label{l.compacidade} $K_{\alpha}\neq\emptyset$ is a compact set.
\end{lemma}

\begin{proof} Let $\{\mu_n\}\subset K_{\alpha}$. By compacity, we can assume that $\mu_n\to
\mu$. Since $V$ is open then $\mu(\Omega\times V)\leq \liminf
(\mu_n)(\Omega\times V)\leq \alpha$. This implies compacity. 
The physical measure given by condition (C2) (see equation (4)) is in
$K_{\alpha}$, because $Leb$-a.e. random orbit stay at most $\gamma_0<\alpha$ inside
$V$ (by $(F1)$). By definition of physical measure (limit of average of Dirac measures
supported on random orbits) and the absolute continuity with respect to the
Lebesgue measure, $\mu_w(V)\leq \alpha$ for $w$ $\BP-$ a.e. holds. In
particular, $\mu(\Omega\times V)\leq\alpha$.  
\end{proof}

We recall that the ergodic decomposition theorem holds for RDS. With this in
mind, we distinguish a set $\SK\subset K_{\alpha}$ :

\begin{defi} $\SK =\{ \mu : \mu_{(w,x)}\in K_{\alpha} \textrm{ for } \mu -a.e.
(w,x) \}$ ($\mu_{(w,x)}$ is the ergodic decomposition of $\mu$).
\end{defi} 

\begin{lemma}\label{l.Kexp} Every measure $\mu\in\SK$ is
$f$-expanding with exponent $c$ :
$$\limsup\limits_{n\rightarrow +\infty}
\frac{1}{n}\sum\limits_{j=0}^{n-1} \log \|Df(T^j(w))(f^j(w)(x))^{-1}\| \leq
-2c $$ for $\mu$-a.e. $(w,x)\in M$. 
\end{lemma}

\begin{proof}We assume first that $\mu$ is ergodic. By definition of $K_{\alpha}$, we have
$\mu(\Omega\times V)\leq\alpha$. But Birkhoff's Ergodic Theorem applied to
$(F,\mu)$ says that in the random orbit of $(w,x)$ $\mu-$a.e. we have: 
$$\lim_{n\to \infty}\frac{1}{n}\sum_{i=0}^{n-1}\chi_V(f^i(w)(x))\leq \alpha.$$ 
Now, we use hypothesis (H1): $\|Df(w,y)^{-1}\|\leq (1+\delta_0)$ for any $y \in
V$ and $\|Df(w,y)^{-1}\|\leq (1+\delta_1)^{-1}$ for any
$y\in V^c$, obtaining:
$$\frac{1}{n}\sum\limits_{j=0}^{n-1} \log \|Df(T^j(w))(f^j(w)(x))^{-1}\| \leq
\log[(1+\delta_0)^{\alpha}(1+\delta_1)^{1-\alpha}] \leq -2c<0$$
$(w,x)-\mu-$a.e.

In the general case we use the ergodic decomposition theorem (see~\cite{O}
and~\cite{LQ}).
\end{proof}

\smallskip

{\it Entropy lemmas}

\begin{defi} Given $\eps >0$, we define :
$$ A_{\eps}(w,x)=\{ y: d(f^n(w)(x), f^n(w)(y))\leq\eps\textrm{ for every }n\geq
0 \} .$$
\end{defi}

\begin{lemma}Suppose that $\mu\in\SK$ is ergodic and let $\eps_0$ given by
lemma~\ref{l.var}. Then, for $\BP$-almost every 
$w$ and any $\eps <\eps_0$,
$$A_{\eps}(w,x)={x} .$$
\end{lemma}

\begin{proof} By lemma~\ref{l.hyptimes} we have infinity hyperbolic times $n_i=n_i(w,x)$ for 
$(w,x)\in H$ (where $\mu(H)=1$). For each $w$ set $H_w=\{x;(w,x)\in H\}$, then
$\BP$-a.e. $w$ we have $\mu_w(H_w)=1$ and infinity hyperbolic times for
$\mu_w$-a.e. $x$. Now, by lemma~\ref{l.var}, if $z\in A_{\eps}(w,x)$ with
$\eps<\eps_0$ we have: 
$$d(x,z)\leq e^{-\frac{cn_i}{2}}d(f^{n_i}(w)(x),f^{n_i}(w)(z))\leq
e^{-\frac{cn_i}{2}}\eps.$$
The lemma follows.
\end{proof}

Let $\SP$ be a partition of $M$ in measurable sets with diameter less than
$\eps_0$. From the above lemma, we get :

\begin{lemma}\label{l.particao} $\SP$ is a generating partition for every $\mu\in\SK$.
\end{lemma}

\begin{proof} As usual we will write:
$$\SP_w^{n}=\{\SC_w^n=(\SP_w)_{i_0}\cap\cdots\cap f^{-(n-1)}(w)(\SP_w)_{i_{n-1}})\}
\textrm{ for each $n\geq 1$},$$
where $(\SP_w)_{i_k}$ is an element of the partition $\SP$.
By the previous lemma, we know that for $\BP$-a.e. $w$, we have
$A_{\eps}(w,x)={x}$ for $x$ $\mu_w$-a.e. Let $A$ a measurable set of $M$ and $\delta>0$. Given $K_1\subset A$ and
$K_2\subset A^c$ two compact sets such that $\mu_w(K_1\bt A)\leq \delta$ and
$\mu_w(K_2\bt A^c)\leq \delta$. Now if $r=d(K_1,K_2)$, the previous lemma says 
that if $n$ is big enough then $diam \SP^n_w(x)\leq \frac{r}{2}$ for $x$ in a
set of $\mu_w$-measure bigger than $1-\delta$. The sets
$(\SC_w^n)_1,\cdots,(\SC_w^n)_k$ that intersects $K_1$ satisfy:
\begin{equation*}\begin{aligned}
\mu( \bigcup (C^n_w)_i\Delta A)
& = \mu(\bigcup (C^n_w)_i-A) + \mu(A - \bigcup (C_w^n)_i)
\\ & \leq \mu(A-K_1)+\mu(A^c-K_2)
  + \delta \leq 3\delta.
\end{aligned}\end{equation*}
This end the proof.

\end{proof}

\begin{corol}\label{c.ks} For every $\mu\in\SK$, $h_{\mu}(f)=h_{\mu}(f,\SP)$
\end{corol}

\begin{proof} The result follows from lemma~\ref{l.particao} and the
theorem~\ref{Kolmogorov-Sinai}.
\end{proof}

We have that the map $\mu\rightarrow h_{\mu}(f,\SP)$ is upper semi-continuous 
at $\mu_0$ measure s.t.
$(\mu_0)_w(\partial P)=0$ for $\BP$-a.e.~$w$, $P\in\SP$. In fact, we have : 
$$h_{\mu}(f,\SP)=\lim\limits_{n\rightarrow\infty}\frac{1}{n}\int
H_{\mu_w}(\SP_w^n) d\BP =
\inf\limits_{n}\frac{1}{n}\int H_{\mu_w}(\SP_w^n) d\BP(w) .$$ But, if
$(\mu_0)_w(\partial P)=0$ for any $P\in\SP$ and $\BP$-a.e. $w$, then the
function $H(\mu, n)$ given by $\mu\rightarrow \int H_{\mu_w}(\SP_w^n) d\BP$ is
upper semi-continuous at $\mu_0$. Indeed, since we are assuming that $T$ is
continuous, the same argument in the proof of
theorem $1.1$ of~\cite{LZ} shows this result. In
particular, because the infimum of a sequence of upper semi-continuous functions
is itself upper semi-continuous, this proves the claim.

\begin{lemma}\label{l.fora} All ergodic measures $\eta$ outside $\SK$ have small entropy :
there exists $\rho_0 <1$ such that
$$h_{\eta}(f)\leq \rho_0 h_{top}(f) .$$ 
\end{lemma}

\begin{proof} By the random versions of Oseledet's theorem and Ruelle's inequality
(see~\cite{Liu}), we have:
$$h_{\eta}(f)\leq \int \sum_{i=1}^{s} \lambda^{(i)}(w,x)m^{(i)}(w,x)d\eta.$$
where $\lambda^{(i)}(w,x)$ and $m^{(i)}(w,x)$ are the Lyapunov exponents of $f$ at
$(w,x)$ and its multiplicity respectively (and
$\lambda^{(1)}(w,x),\cdots,\lambda^{(s)}(w,x)$ are the positive Lyapunov
exponents). 
Furthermore, by hypothesis the measure is ergodic, then these objects are
constant a.e. then $h_{\eta}(f)\leq \sum_{i=1}^s \lambda^{(i)}$ and $\int \log
\|\det Df(w)(x)\|d\eta=\sum_i \lambda^{(i)}$. Since $\|Df(w)(x)^{-1}\|\leq
(1+\delta_0)$ we have $\lambda_l>-\log (1+\delta_0)$. By the definitions of $m_2$,
$M_2$ and the above estimates, we have by (C1):
\begin{eqnarray*}
h_{\eta}(f) &\leq & \int \log \|Df(w)(x)\|d\eta - \sum_{i=s+1}^l \lambda^{i} \\
& \leq & \eta(\Omega\times V)m_2 +(1-\eta(\Omega\times V))M_2 +(l-s)(1+\delta_0)\\
& \leq & \alpha m_2 + (1-\alpha)M_2 + l\log(1+\delta_0)
\end{eqnarray*}
Now the physical measure $\mu_{\BP}$ given by condition (C2) satisfy
$\mu_{\BP}(W)<\gamma_0$ (by $(F1)$). The Random Pesin's formulae gives:
$$h_{\mu_{\BP}}(f)=\int\log\|det Df\|d\mu_{\BP}\geq
\mu_{\BP}(W)m_1+(1-\mu_{\BP}(W))M_1.$$
But $m_1<M_1$ then $\gamma_0m_1+(1-\gamma_0)M_1\leq h_{\mu_{\BP}}(f)$. 
Using that $\eta \notin K$, $m_2<M_2$ and (C1) we have:
$$\alpha m_2 + (1-\alpha)M_2 < \gamma_0 m_1 + (1-\gamma_0)M_1 -
l\log(1+\delta_0).$$
Then, we can choose $\rho_0<1$ such that
$$\alpha m_2 + (1-\alpha)M_2 + l\log(1+\delta_0)< \rho_0(\gamma_0 m_1 +
(1-\gamma_0)M_1 )< \rho_0 h_{\mu_{\BP}}(f)$$
This gives: $h_\eta(f) \leq \rho_0 h_{top}(f)$.
\end{proof}

\begin{corol} $\pi_F(\phi)= \sup\limits_{\mu\in\SK}\{ h_{\mu}(f)+\int\phi
d\mu\}.$
\end{corol}

\begin{proof} By remark~\ref{r.erg}, we need to show that:
$$\sup_{\mu\in\SK} \{h_{\mu}(f)+\int\phi d\mu\}=\sup_{\mu\in\SM_e(\mundo)}
\{h_{\mu}(f) +\int\phi d\mu\}$$
By the previous lemma, if $\eta \notin \SK$ then:
$$h_{\eta}(f)+\int \phi d\eta \leq \rho_0h_{top}(f)+\|\phi\|_1<\pi_f(\phi)$$
\end{proof}

\begin{proof}[Proof of theorem A] We will use the following notation: $\Psi(\mu)=h_{\mu}(f)+\int\phi d\mu$.
Let $\{\mu_k\}\subset \SK$ such that $\Psi(\mu_k)\to\pi_f(\phi)$, by compacity
we can suppose that $\mu_k$ converge to $\mu$ weakly. 

Fix $\SP$ a partition with diameter less than $\eps_0$, and for $w$-a.e., $\mu_w(
\partial P)=0$, for any $P\in\SP$. By corollary~\ref{c.ks} we have $h_{\mu_k}(f)=h_{\mu_k}(f,\SP)$.
Then $\pi_f(\phi)=\sup\limits_{\eta \in \SK}\Psi(\eta)=\limsup \Psi(\mu_k)$. By
the comments after corollary~\ref{c.ks} we know that $\eta \to h_{\eta}(f,\SP)$
is upper semicontinuous in $\eta$ {\em over $\SK$}, then:
$$\limsup \Psi(\mu_k) \leq h_\mu(f,\SP)+\int\phi d\mu\leq \Psi(\mu).$$
But, $\Psi(\mu)\leq \pi_f(\phi)$. This implies that $\mu$ is an equilibrium
state.

In the other hand, if $\eta$ is a measure which attain the supremum
in~(\ref{e.pressure}) then let $\eta_{(w,x)}$ the ergodic decomposition of $\eta$.
Then the entropy of $\eta$ is equal to the integral of entropies of its ergodic
components (see~\cite{Liu}, page 1289 and references there in), of course the same occurs
with the $\Psi(\eta)$ (*). If $(x,w)\notin\{(x,w);\eta_{(x,w)}\in \SK_{\alpha}\}$ then
by lemma~\ref{l.fora}:
$$\Psi(\eta_x)=h_{\eta_{(x,w)}}(f)+\int\phi d\eta_{(x,w)}\leq
\rho_0h_{top}(f)+\|\phi\|_1<\pi_f(\phi).$$
Then if $\eta(\{(x,w);\eta_{(x,w)}\in \SK_{\alpha}\}^c)>0$, (*) says that
$\Psi(\eta)<\pi_f(\eta)$ a contradiction, so every equilibrium state is in
$\SK$. The proof of the theorem is now complete.
\end{proof}

\begin{rem} Liu-Zhao~\cite{LZ} show the semi-continuity of the entropy under
the hypothesis that $T:\Omega\rightarrow\Omega$ is \emph{continuous} and $f$ is
expansive \emph{at very point of $M$}. From this
result, a natural
question is : ``What about the semi-continuity without topological assumptions
(e.g., continuity) ? And the case of weak expansiveness assumptions ?''. We 
point
out that the proof of theorem $A$ shows the semicontinuity of the entropy map 
\emph{in the set
$\SK$}. This partially answer the question since, although we need to assume
continuity, only a weak expansion at Lebesgue a.e. point of $M$ is required 
(\emph{this assumption is the sole reason of the restriction to the set of 
measures $\SK$}). Indeed, non-uniform expansion on Lebesgue a.e. point obligates us to
restrict the proof of our lemmas on semicontinuity to the set $\SK$.  
\end{rem}

\begin{rem}Our theorem $A$ holds in the context of \emph{RDS bundles} (see~\cite{Liu}
or~\cite{LZ}) with the extra assumption that $T$ \emph{and} the skew-product $F$ are
continuous.
\end{rem}

\section{Appendix}
\label{Apendice}
 
We now prove that ($F1$) follows from (H1) and (H2), in fact, this is a well know
argument (see for example~\cite{AA}), but for sake of completeness we give the proof.

Fix $(w,x)$, if $i=(i_0,\cdots i_{n-1})\in \{1,\cdots,p+q\}^n$ let
$[i]=B_{i_0}\cap f^{-1}(w)(B_{i_1})\cap \cdots\cap f^{-n+1}(w)(B_{i_{n-1}})$ and
$g(i)=\#\{0\leq j<n; I_j\leq p\}$.

If $\gamma>0$ then $\#\{i;g(i)<\gamma n\}\leq \sum\limits_{k\leq\gamma n}{n
\choose k}p^{\gamma n}q^n.$ By Stirling's formula this is bounded by
$(e^{\xi}p^{\gamma}q)^n$ (here $\xi$ depends of $\gamma$) and $\xi(\gamma)\to 0$
if $\gamma \to 0$.

Now (H1) and (H2) says that $m([i])\leq \sigma_1^{-n}\sigma_1^{\gamma n}$. If we
set $I(n,w)=\bigcup \{[i]; g([i])<\gamma n\}$ then $m(I(n,w))\leq
\sigma_1^{-(1-\gamma)n}(e^{\xi}p^{\gamma}q)^n$ and since $\sigma_1>q$ there is a
$\gamma_0$ (small) such that $(e^{\xi}p^{\gamma}q)^n<\sigma_1^{(1-\gamma_0)}$.
Then there is a $\tau=\tau(\gamma_0)<1$ and $N=N(\gamma_0)$ such that if $n\geq
N$ then $m(I(n,w))\leq \tau^n$

Let $I_n=\bigcup\limits_{w}(\{w\}\times I(n,w))$ and by Fubini's theorem
$\BP\times Leb (I_n)\leq \tau^n$ if $n\geq N$. But $\sum\limits_n \BP\times Leb
(I_n)<\infty$ then Borel-Cantelli's lemma implies:
$$\BP\times Leb(\bigcap\limits_{n\geq 1}\bigcup\limits_{n\geq k} I_k)=0$$
Using Birkhoff's theorem we have that the set:
\begin{eqnarray*}
\{(w,x);\exists \ n\geq 1, \forall k\geq n, \lim \frac{\#\{0\leq j<n;
f^j(w)(x)\in B_1\cup\dots\cup B_p\}}{n}\}.
\end{eqnarray*}
has $\BP\times Leb$-measure at least $\gamma_0$. Now by Fubini's theorem again, we have (F1).

\textbf{Acknowledgements.} The authors are indebted to Professor Marcelo Viana 
for useful conversations, suggestions and advices. Also to
Professor Pei-Dong Liu for give
us helpful references (including the preprint~\cite{LZ}), for communicated to us 
some entropy results and for some corrections on old versions of the paper, and Professor Yuri
Kifer for another valuable references. Finally, we are grateful to IMPA and his
staff for the remarkable research ambient.

\smallskip


\vfill

{\footnotesize
\noindent 

\vfill

Alexander Arbieto ({\tt alexande{\@@}impa.br}) \\

Carlos Matheus ({\tt cmateus{\@@}impa.br}) \\

Krerley Oliveira ({\tt krerley{\@@}mat.ufal.br}) \\
\smallskip

\noindent IMPA, Estrada D. Castorina 110, Jardim Bot\^anico, 22460-320 Rio de Janeiro, Brazil
}

\end{document}